\documentclass[a4paper,10pt]{article}

\newtheorem{thm}{Theorem}
\newtheorem{cor}[thm]{Corollary}
\newtheorem{rem}[thm]{Remark}
\newtheorem{lem}[thm]{Lemma}
\newtheorem{prop}[thm]{Proposition}
\def \a{{\alpha}}
\def \b{{\beta}}

\def \DD{{\mathcal D}}
\def \d{{\delta}}
\def \e{{\varepsilon}}

\def \g{{\gamma}}
\def \G{{\Gamma}}

\def \l{{\lambda}}

\def \o{{\omega}}
\def \O{{\Omega}}
\def \p{{\varphi}}
\def \t{{\vartheta}}

\def \m{{\mu}}
\def \s{{\sigma}}

\def \A{{\cal A}}

\def \E{{\bf E}\, }

\def \N{{\bf N}}

\def \P{{\bf P}}

\def \qq{{\qquad}}
\def \R{{\bf R}}

\def \T{{\bf T}}
\def \td{{\widetilde D}}
\def \ua{{\underline{a}}}
\def \ue{{\underline{\e}}}

\def \uz{{\underline{z}}}

\def \dd{{\rm d}} 
\def \noi{{\noindent}}
\def \nh{{\hat \nu}} 
\def\qed{\hbox{\vrule height 6pt depth 0pt width
6pt}}
\def\cqfd{\hfill\penalty 500\kern 10pt\qed\medbreak}
 \font\phh=cmcsc10 at  8 pt
 \scrollmode
  \font\gsec= cmb10 at 10 pt 

 \scrollmode

  at 11,6 pt

   \title{ \bf Dirichlet polynomials: some old and recent results, and  their interplay in number theory}
 \author{
   Michel WEBER}

\begin{document}
 
\maketitle

\begin{abstract}
 In the first part of this expository paper, we present and discuss the interplay of  Dirichlet polynomials in some classical problems
of number theory, notably the Lindel\"of Hypothesis. We review some typical   properties of their means and continue with some
investigations concerning their supremum properties. Their random counterpart is   considered in the last part of the paper, where a
analysis of their supremum properties,  based on methods of stochastic processes, is developed.  
\end{abstract}


\section{Introduction}
 This is an expository   article on the interplay of  Dirichlet polynomials in some classical problems
of number theory, notably the Lindel\"of Hypothesis (LH), the Riemann Hypothesis (RH), as well as    on their typical means and   supremum
properties.   
Some of the efficient methods  used in this context are also sketched. In  recent works
\cite{LW1},\cite{LW2},\cite{LW3},\cite{W4},\cite{W5},\cite{W6},\cite{W7}    lying at the interface of probability theory, the
theory of Dirichlet polynomials and of the one of the   Riemann zeta-function, we had to combine with these questions. The
growing interaction between various specialities of the analysis,  further  motivated us  in this project to put in the same framework a
certain number of basic and very important results and tools arising from the theory of Dirichlet polynomials and of the    Riemann
zeta-function; and to let them at disposal to   analysts and probabilists who are not necessarily number theorists. In doing so, our
wish is to spare their time in the sometimes  tedious enterprise of finding the relevant results with the mostly appropriate methods to
establish them.

The interplay with the LH and RH is presented in Section 2, where some
equivalent reformulations of the LH in terms of approximating   Dirichlet polynomials,   arising notably from   T\'uran's works, are
discussed. The link between the RH and the absence of zeros  of the   approximating  
Dirichlet polynomials  in some regions of the complex plane was thoroughly investigated  by T\'uran
in \cite{Turan},\cite{Turan1},\cite{Turan2},\cite{Turan3}  and later by Montgomery  in
\cite{[Mon]},\cite{[Mon1]}.  Some of the  most striking results are presented. 

In Section 3, we investigate the behavior of the mean value of Dirichlet polynomials. The used reference sources are  \cite{Iv},
\cite{[Mon1]},  \cite{Ra} and naturally \cite{Ti}. To begin, we follow an approach based in the Fourier inversion formula. Next, the
mean value estimates  are established by means of a version of Hilbert's inequality due to Montgomery and Vaughan \cite{[Mon]}. A simple
argument is provided  for establishing the  lower bound.   The corresponding results  for the zeta-function are briefly mentionned and
discussed.    Some basic results concerning suprema of Dirichlet polynomials are presented in Section 4.  The analysis of the suprema of
their random counterpart is made in the two last sections. These are notably built on works of Hal\'asz \cite{Ha1}, Bayart, Konyagin,
Queff\'elec
\cite{BKQ},\cite{KQ},\cite{Q1},\cite{Q2},\cite{Q3} and Lisfhits, Weber
 \cite{LW1},\cite{LW2},\cite{W4}. 
\section{Interplay in Number Theory}
 To any real
valued function
$d$    defined on the integers, we may associate the Dirichlet polynomials 
\begin{equation} \label{e12}
D_N(s)= \sum_{n=1}^N  {d(n)\over  n^{ s}},   \qq\qq\qq (s=\s+it).
\end{equation}  
The particular case $d(n)\equiv 1$ is already of crucial importance, since it is intimately related   to the behavior of the Riemann
Zeta-function. We shall first investigate this link. Recall that the Riemann
Zeta-function is  defined
on the half-plane $\{s: \Re s>1\}$ by the series
 $\zeta(s)= \sum_{n=1}^\infty n^{-s}$,  
which admits a meromorphic continuation to the entire complex plane. 
\smallskip\par
And we  have  the following classical approximation result
(\cite{Ti}, Theorem 4.11)   
\begin{equation}\label{approx} \zeta(s) =\sum_{n\le x} {1\over n^s} -{x^{1-s}\over 1-s} +{\cal O}(
x^{-\s}),
\end{equation}
uniformly in the region $\s\ge \s_0>0$, $|t|\le T_x:=2\pi x/C$, $C$ being a constant
$>1$. The celebrated 
Lindel\"of Hypothesis   claiming that
\begin{equation}\zeta( {1\over 2}+it) ={\cal O}_\e (t^\e)   
\end{equation}
 can be reformulated in terms of Dirichlet polynomials. This was observed since quite a long time and  Tur\'an \cite{Turan1} had  shown
that the truth of the inequality 
\begin{equation}\label{alt}\Big|\sum_{n=1}^N {(-1)^n \over n^{it}} \Big|\le C N^{1/2+\e}(2+|t|)^ \e
 \end{equation}
with an arbitrary small $\e>0$, is in turn  equivalent to the  LH. 
  Alternatively, the  equivalent reformulation  (\cite{Ti},
Chap. XIII)
\begin{equation}\label{equivdemi}{1\over T}\int_1^T \big|\zeta({1\over 2} + it)\big|^{2k} {\rm d}t= {\cal
O}_\e\big( T^\e\big), \quad\qq k=1,2,\ldots 
\end{equation}
  which reduces to  
\begin{equation} \int_{T/ 2}^T \big|\zeta({1\over 2} + it)\big|^{2k} {\rm d}t={\cal O}_\e\big(
T^{1+\e}\big),
\quad\qq k=1,2,\ldots
 \end{equation}
 yields  when combined with (\ref{approx}), the equivalence with
\begin{equation} \label{reducLH}\int_{ n/2}^{ n} \Big| \sum_{m= 1}^{n}  {1\over  m^{1/2+it }}\Big|^{2k} {\rm d}t= {\cal
O}\big(  n^{1+\e}\big), \quad\qq k=1,2,\ldots \end{equation}
 or, by  using  Euler-MacLaurin formula,  with  
\begin{equation} \int_{ n/2}^{ n} \Big| s
\int_{  0}^{n}   y^{-1-s }B_1(y) \dd y\Big|^{2k} {\rm d}t= {\cal
O}\big(  n^{1+\e}\big), \quad\qq k=1,2,\ldots
\end{equation}
  where $B_1(y)=\{ y\} -{1/ 2}$, mod $1$ is the first Bernoulli function. 
   The term
${n^{1-s}/( 1-s)}$ can be indeed  neglected
 since 
$$ \int_{n\over 2}^n\Big|
 {n^{{1\over 2}-  it}\over {1\over 2}-  it }\Big |^{2k}{\rm d}  t \le Cn^{k}\int_{n\over 2}^{\infty}      {\dd t  \over
  {({1\over 4}+t^2)^{k}}} 
\le   C_k\,n^{k }n^{ -2k+1}\le   C_k\, n^{    1- k } . $$

We see with (\ref{reducLH}) that all the mystery of the LH is
hidden in the Dirichlet sum $\sum_{m= 1}^{n}     m^{ -1/2-it } $, $n/2\le t\le n$.

\begin{rem} \rm The
best known result   is due to Huxley \cite{[H]},
   \begin{equation}\zeta( {1\over 2}+it) ={\cal O}_\e(t^{32/205 +\e}) .
   \end{equation}
In the recent work \cite{LW3}, in order to understand the behavior of $\zeta({1\over 2} + it)$ as $t $ tends to infinity,    
the time
$t$ is modelized  by a  Cauchy random walk, namely the sequence of partial sums $S_n= X_1+ \ldots +X_n   
$ of a   
sequence of independent Cauchy distributed random variables $X_1,X_2,\ldots$ (with characteristic function $\varphi(t) =e^{-|t|}$).  The
almost sure asymptotic behavior of the system  
$$  \zeta_n:=\zeta({1\over 2} + iS_n), \quad\qq n=1,2, \ldots  
$$
is  investigated. 
   Put for any positive integer $n$ 
$${\mathcal Z}_n = {\zeta}({1/  2}+iS_n) -\E\,  {\zeta}({1/ 2}+iS_n)
=\zeta_n-\E \zeta_n.  
$$
The  crucial  preliminary study of  second order properties of the system $\{ {\mathcal Z}_n,n\ge 1\}$  
yields the  
   striking fact   that  this one nearly behaves like a system of non-correlated 
variables, i.e. the variables ${\mathcal Z}_n$ are weakly orthogonal. More precisely,  there exist constants $C,C_0$  
\begin{eqnarray}\E\,   | {\mathcal Z}_n|^2
&=&\log n + C   
 + o(1), \qquad n\to\infty,
\cr   
\hbox{and for $m>n+1$, \qq}\big|\E\,{\mathcal Z}_n\overline{{\mathcal Z}_m}\big|
&  \le &C_0 \max \Big( {1\over n},{1\over 2^{ m-n  }}   \Big). \qq\qq\qq\qq
 \end{eqnarray}
 The proof is very technical. 
And the main result of \cite{LW3}, which follows   from a convergence criterion,  states
  \begin{equation}\lim_{n\to \infty}{\sum_{k=1}^n \zeta({1\over 2} + iS_k) - n  \over n^{{1/2}}(\log n)^b
}\buildrel{(a.s.)}\over=0 ,
\end{equation}
for any real $b>2$.
\end{rem}
\bigskip\par 
 Now we pass to the interplay with  the Riemann Hypothesis. 
As is well-known, the Riemann Zeta-function has a unique and simple pole  of residue $1$ at $s=1$, and simple zeros at
$-2,-4,-6,\ldots$, and only at these points which are called trivial zeros.
There exist also non-trivial zeros in the band  $\{s: 0< \Re s < 1\}$.
   The
Riemann Hypothesis
  asserts that all non-trivial zeros of the   $\zeta$-function have
abscissa ${1\over 2}$.   The validity of
RH  implies (\cite{Ti},   p.300) that
$$\zeta( {1\over 2}+it) ={\cal O}\Big(\exp\big\{ A{\log t\over \log\log
t}\big\}\Big), 
$$
$A$ being a constant, which is even a stronger form of
LH; the latter  being strictly weaker than RH.
 
In several papers, Tur\'an   investigated the interconnection between the RH and the absence of zeros  of the   approximating  
Dirichlet polynomials  in some regions of the complex plane. For instance, he proved in   \cite{Turan}  that if for   $n>n_0$,
none of the Dirichlet polynomials $\sum_{m=1}^n m^{-s}$ vanishes in a half-strip
$$\s\ge 1+ {\log^3 n \over \sqrt n}, \qq \gamma_n\le t\le \gamma_n + e^{n^3},$$
with a suitable real $\gamma_n$, then the RH is true.  

However  Montgomery  has shown in \cite{[Mon1]} that for $n>n_0$, any interval $[\gamma, \gamma  + e^{n^3}]$ contains the imaginary part
of a zero.
  Another definitive result established by Montgomery in the same paper is that if $0<c<{4\over \pi} -1$, then for all $n>n(c)$,
$\sum_{m=1}^n m^{-s}$ has zeros in the half-plane
\begin{equation}\label{zfr2}\s >1+c{\log\log n\over \log n}, 
\end{equation}
whereas if $ c>{4\over \pi} -1$, $n>n(c)$, then $\big|\sum_{m=1}^n m^{-s}- \zeta(s)\big|\le  |\zeta(s)|/2 $, so that
$\sum_{m=1}^n m^{-s}$ do not vanish in a half-strip (\ref{zfr2}).
\smallskip\par

In the other direction, Tur\'an showed  (\cite{Turan3}, Satz II) that if the RH is true, none of the Dirichlet
polynomials $\sum_{m=1}^n m^{-s}$ vanishes in a half-strip
\begin{equation}\s\ge 1 , \qq c_1\le  t\le e^{c_2\sqrt{\log n \log\log n}},\end{equation}when $n\ge c_3$.

 There are also   results for Dirichlet
polynomials expanded over the primes. Their local suprema are intimately connected with zerofree regions of the Riemann
zeta-function. Among the several results proved in \cite{Turan2}, we may quote the following. Suppose there are constants
$\a\ge 2$,
$0<
\b\le 1$,
$
\tau(\a,
\b)$, such that for a
$\tau >
\tau(\a,
\b)$ the inequality
 \begin{equation} \label{zfr}\Big|\sum_{N_1\le p\le N_2 }  {1\over p^{i\tau}}\Big|\le {N\log^{10} N\over \tau^\b}
 \end{equation}   
  holds for all $N_1,N_2$ integers with
 $\tau^\a\le N\le N_1\le N_2\le 2N\le e^{\tau^{\b/10}}  $. Then (\cite{Turan2}, Theorem 1) $\zeta(\s +i\tau)$ does not vanish if
$\s>1-\b^3/(e^{10}\a^2)$. 
For the sake of orientation, Tur\'an also remarked that for the sum 
$$S=\sum_{N\le n\le 2N }  {1 \over n^{it}}, \qq \tau \ge 2 $$
the elementary formula
 $|(\nu+1)^{1+i\tau}-\nu^{1+i\tau}- ( {1+i\tau})\nu^{ i\tau}|\ll \tau^2/\nu  $,
gives at once
\begin{equation}  |S|\ll (N\log N)/\tau
 \end{equation}
if $N\ge \tau^2$.  But the relevant sum is  $\sum_{n=1}^N {(-1)^n \over n^{it}}$ according to (\ref{alt}).
       
 Farag recently showed in \cite{Fa} that these sums possess zeros near every vertical line in the critical strip. The proof is notably based on a
"localized" version of Kronecker's Theorem (section 3).
\smallskip\par
 The likely best known result concerning zerofree regions is due to Ford \cite{Fo}: $\zeta(\sigma + it)  \not= 0$ whenever $|t| \ge 3$ and
 \begin{equation}    \sigma\ge 1-\frac{1}{57.54(\log{|t|})^{2/3}(\log{\log{|t|}})^{1/3}}. 
  \end{equation}
\begin{rem}\rm Speaking of the RH,  it is difficult not mentionning the striking equivalent
reformulation proved by  Robin in
\cite{Ro}, which is  
  at the same time likely   the most simple.   Let an integer $n$ be termed ''colossally abundant'' if, for
some
$\e > 0$,
$\s(n)/n^{ 1+\e}\ge 
\s(m)/m^{ 1+\e}$ for
$m < n$ and $\s(n)/n^{ 1+\e}> 
\s(m)/m^{ 1+\e}$ for $m > n$,   where $\s(n)$ is the sum of divisors of $n$. Using colossally abundant numbers,  Robin showed that the 
RH  is true if and only if 
$$ {\s(n)\over n} < e^\gamma\log\log n,$$
for $n>5040$, where   $\gamma$ is Euler's constant. 
Let  $\{x_n, n\ge 1\}$ be the sequence
of colossally abundant numbers. In the same paper, he also showed that the sequence
$\{\s(x_n)/x_n \log \log  x_n, n\ge 1\}$ contains an infinite number of local extrema.   
In relation with Robin's result,  Lagarias  
showed in \cite{La} that the RH  is true if and only if
$$ \s(n) \le  H_n +e^{H_n} \log H_n,$$
where $H_n =
\sum_{
j\le n}  1/j$ is the
$n$-th  harmonic number. 

\smallskip\par 
  
  Grytczuk \cite{Gr} investigated the upper bound for $\s(n)$ with some different
$n$. Let $(2, n) = 1$ and $n =
\prod_{j=1}^k
 p_j^{\a_j}$ , where the $p_j$ are
prime numbers and $\a_j\ge 1$. Then, for all odd positive integers $n > 3^9 /2$,
$$\s(2n) < {39\over
40} e^{\gamma} 2n\log \log 2n,\qq \hbox{and}\qq \s( n) <   e^{\gamma}  n\log \log  n.$$
  Some other criteria equivalent to the RH can be found in  \cite{CW}.
\end{rem}

 \section{Mean Values of Dirichlet Polynomials}
  Let  $k$ be some positive integer. Considerable efforts  were made to finding good estimates  for   the mean integrals 
$${1\over 2T}\int_{-T}^T \Big|  \sum_{n=1}^N  {1\over  n^{ \s+it }}\Big|^{2k}\dd t,  $$
because of  Bohr's theory of
almost periodic functions, and of the "mean value" equivalent   reformulations of the LH.  First notice that
\begin{equation}\label{lim}\lim_{T\to\infty}{1\over 2T}\int_{-T}^T \Big|  \sum_{n=1}^N  {1\over  n^{ \s+it }}\Big|^{2k}\dd t
=\sum_{m=1}^\infty {d_{k,N }^2(m)\over m^{2\s}}, 
\end{equation} 
 where $d_{k,N }(m)$ denotes the number of representations of  $m$ as a product of
$k$ factors less or equal to $N$.  We propose to deduce this 
from the Fourier inversion formula. This seems a natural approach, although we could not refer to some book or paper.
 If  
$\nu$ is a distribution function on $\R$ and   $\widehat  \nu (t)=\int_\R e^{it x} \nu(dx) $ denotes its characteristic function, then
\begin{equation}\label{fif0}\lim_{T\to \infty}{1\over 2T}\int_{-T}^Te^{-itx_0}\widehat\nu(t) dt=\nu\{x_0\}. 
\end{equation}
From this  also follows that 
 $\lim_{T\to \infty} {1\over 2T}\int_{-T}^T |\widehat \nu(t)|^2dt=\sum_{x\in
\R}\nu(\{x\})^2$ and more generally    
\begin{equation}\label{fif}\lim_{T\to \infty} {1\over 2T}\int_{-T}^T |\widehat \nu(t)|^{2k}dt =\sum_{x\in
\R} \nu  ^{*k}(\{x\})^2,
\end{equation}  
for 	any    positive integer $k $.  Apply (\ref{fif}) to the measure
    $\nu =\sum_{n=1}^N {1\over n^\s}  
\delta_{\{-\log n\}},  $   where $\delta_{\{x\}}$ is the  Dirac measure at point  $x$, then 
 $ \widehat \nu(t) =
\sum_{n=1}^N   n^{-(\s+ it)}$ and
\begin{eqnarray*}& &\!\!\!\!\!\!\!\!\!\!\!\! \lim_{T\to \infty} {1\over 2T}\int_{-T}^T \big|
\sum_{n=1}^N {1\over n^{\s+ it}} \big|^{2k}dt    =  \sum_{x\in \R}\Big(\sum_{   n_1 \ldots n_k  =e^x} {1\over n_1^{\s }\ldots n_k^{\s
}}\Big)^2 \cr &=&
 \sum_{    Y= e^x\atop Y\in \N}   {\#\{ (n_1,\ldots, n_k), n_i\le N:  \prod_{i=1}^kn_i=Y    \}^2
\over  Y^{ 2\s }}   
  = 
  \sum_{    Y \in \N}{d_{k, N }^2(Y)\over  Y^{ 2\s }}.  
\end{eqnarray*} 

 Clearly
$$\lim_{N\to \infty} \sum_{    Y \in \N}{d_{k, N }^2(Y)\over  Y^{ 2\s }}=\sum_{m=1}^\infty {d_k^2(m)\over m^{2\s}}.  $$

As an equivalent formulation of the LH is that for $\s>{1\over 2}$, $ k=1,2,\ldots$
\begin{equation}\label{LH2}\lim_{T\to\infty}
{1\over T}\int_1^T \big|\zeta(\s + it)\big|^{2k} {\rm d}t
 =
\sum_{m=1}^\infty {d_k^2(m)\over m^{2\s}} 
\end{equation}
where $d_k(m)$ denotes the number of representations of   $m$ as a product of
$k$ factors, the LH can be interpreted as a kind of generalized Fourier inversion formula for the infinite measure $ \sum_{n=1}^\infty
{1\over n^\s}  
\delta_{\{-\log n\}}   $.  This approach is investigated in \cite{W5}, see also \cite{W8} section 13.7.
\smallskip\par  
One  can precise    a little  
 more   (\ref{fif0}). Let 
$$ {\bf M}_T(\nu ,x_0) ={1\over 2T}\int_{-T}^Te^{-itx_0}\widehat\nu(t) dt .$$

\begin{prop} For  any    
arbitrary non-decreasing sequence
$ \{T_p, p\ge 1\}$    of positive reals,  any
$x_0$    
\begin{equation}\label{sqfif}\sum_{k=1}^\infty \big|{\bf M}_{T_{k+1}} (\nu,x_0) -{\bf M}_{T_{k }} (\nu,x_0) \big|^2\le 24 \nu^2 (\R).
 \end{equation}
\end{prop} 
 The total mass of the measure appears this time, unlike in (\ref{fif0}).
Notice also that (\ref{sqfif}) alone    already implies  that ${\bf M}_T(\nu ,x_0)$    converges,    as $T$ tends to
infinity. 
\medskip\par\noi   
   {\it Proof.}  Consider the kernels  
 $${\cal V}_T(\t)={e^{iT\t}-1\over iT\t},\qq\quad V_T(y)=\Re \big({\cal V}_T(y)\big)={\sin Ty\over Ty}. $$
 By the Cauchy-Schwarz inequality, we first observe that
\begin{eqnarray*} & &\big|{\bf M}_{T_2}(x_0)-{\bf M}_{T_1}(x_0)\big|^2 = \Big|\int_\R
\big[V_{T_2}(x-x_0) -V_{T_1}(x-x_0)
\big]\nu(dx)\Big|^2\cr &\le& \nu(\R).\int_\R \big[V_{T_2}(x-x_0) -V_{T_1}(x-x_0)
\big]^2\nu(dx).
\end{eqnarray*}
Introduce a new measure $\nh$, a regularization of $\nu$ defined as follows:
$$ {{d\nh}\over {dx}} (x) =
   \int_{\vert \t\vert<\vert x \vert} \vert x \vert^{-3}\t^2
\nu(d\t)
  +\int_{\vert \t \vert\ge \vert x \vert   } \vert \t \vert^{-1}
\nu(d\t) . $$
 From the basic elementary inequalities
\begin{eqnarray*}|{\cal V}_{T_2}(\t)-{\cal V}_{T_1}(\t)|& \le &\min \Big\{
{T_2-T_1 \over 2} |\t| , {2(T_2-T_1)\over T_2} \Big\},\qquad T_2\ge T_1,\cr   
  |{\cal V}_{T_1}(\t)|&\le &{2\over T_1|\t|},
\end{eqnarray*}
   we   have (\cite{LW4}, Section 4)
$$\Vert {\cal V}_{T_2}-{\cal V}_{T_1}\Vert_{2,\nu}^2 \le
8\nh({1\over T_2},{1\over T_1}]. 
$$
 Thereby
 \begin{eqnarray*}  \int_\R \big[V_{T_2}(x-x_0) -V_{T_1}(x-x_0)
\big]^2\nu(dx)   &\le &  \int_\R \big|{\cal V}_{T_2}(y) -{\cal V}_{T_1}(y)
\big|^2\nu_{x_0}(dy) \cr&\le& 8\widehat \nu_{x_0} ({1\over T_2},{1\over T_1}],   
\end{eqnarray*}
 where   we write  $\nu_{y}(A)=\nu (A-y) $, for each  $A\in \mathcal B({\bf R})$. 
And
$$ \widehat \nu_{x_0}(\R)    =  
  \int \left(
\int_{\vert \t\vert<\vert x \vert} \vert x \vert^{-3} dx\  \t^2
  +\int_{\vert x \vert<\vert \t \vert }  dx \vert \t \vert^{-1}
\right) \nu(d\t) \le
 3   \nu_{x_0}(\R) = 3  \nu (\R).  $$
Therefore
$$\big|{\bf M}_{T_2}(x_0)-{\bf M}_{T_1}(x_0)\big|^2\le 24\nh({1\over T_2},{1\over T_1}]. $$
The claimed inequality follows easily.
 \medskip\par \cqfd 

    Let    $k,N$ be fixed but arbitrary positive integers. Put for $T>0$
$$  M _{T } = {1\over 2T}\int_{-T}^T \Big|  \sum_{n=1}^N  {1\over  n^{ \s+it }}\Big|^{2k}\dd t.$$
As an immediate consequence of the preceding Proposition, we get 
\begin{cor}  
For  any    
arbitrary non-decreasing sequence
$ \{T_j, j\ge 1\}$    of positive reals, 
\begin{equation}\label{sqfif}\sum_{j=1}^\infty \big|M_{T_{j+1}}   -M_{T_{j }}    \big|^2\le 3.2^{2k+3}\Big( \sum_{n=1}^N{1\over
n^{\s}}\Big)^{2k} .
 \end{equation} 

\end{cor}

\bigskip\par  

Now let 
$ \l_1, \ldots, \l_N$ be distinct  real numbers and consider the simplest mean integral 
$$ \int_0^T \Big|  \sum_{n=1}^N   d(n)   e^{  it\l_n } \Big|^2\dd t.$$
  The following form of Hilbert's inequality due to Montgomery and Vaughan yields   precise  estimates  of this integral. 
\begin{lem} \label{HMV} Let 
$\d >0$  be   real number  such that   $|\l_m-\l_n|\ge \d $ whenever $m\not=n$. Then 
\begin{equation}\label{HI} \bigg|\sum_{1\le m,n\le N\atop m\not=n}{x_my_n\over \l_m-\l_n}\bigg|\le {\pi\over
\d}\Big(\sum_{m=1}^N|x_m|^2\Big)^{1/2}\Big(\sum_{n=1}^N|y_n|^2\Big)^{1/2}.
\end{equation}
\end{lem} 
By  squaring out and integrating term-by-term, we get
$$ \int_0^T \Big|  \sum_{n=1}^N   d(n)   e^{  it\l_n } \Big|^2\dd t=  T\sum_{n=1}^N  d^2(n) +2\Re\Big\{
\sum_{1\le m<n\le N}d(m)\overline{d(n)} \,
{e^{i(\l_m-\l_n)T}-1\over i(\l_m-\l_n)}\Big\}.$$
Since the sum in the parenthesis is the difference 
$$ \sum_{1\le m<n\le N}d(m)e^{i \l_m T}\overline{d(n)} \,e^{-i  \l_n T} 
{1 \over i(\l_m-\l_n)}-\sum_{1\le m<n\le N} 
{ d(m)\overline{d(n)} \, \over i(\l_m-\l_n)},$$
  by applying (\ref{HI}) to each part, we get
$$ \Big|
\sum_{1\le m<n\le N}d(m)\overline{d(n)} 
{e^{i(\l_m-\l_n)T}-1\over i(\l_m-\l_n)}\Big|\le  {2\pi\over
\d} \sum_{n=1}^N|d(n)|^2   .$$

Consequently
\begin{prop}\label{MV}\begin{equation}   \bigg|\int_0^T \Big|  \sum_{n=1}^N   d(n)   e^{  it\l_n } \Big|^2\dd t -  T\sum_{n=1}^N 
|d(n)|^2 \bigg|\le   {4\pi\over
\d} \sum_{n=1}^N|d(n)|^2  .
\end{equation}
\end{prop}
 
Choose   $\l_n=\log n$ and observe that for $m<n\le N$, $\l_n -\l_m= \log {n\over m}\ge \log {n\over n-1 } \ge c N^{-1}$. 
We get in this case\begin{equation}\label{MVD}  \bigg|{1\over T}\int_0^T \Big|  \sum_{n=1}^N   {d(n)\over   n^{  it  }} 
\Big|^2\dd t -   \sum_{n=1}^N 
d^2(n)
\bigg|\le  C{ N\over T}   \sum_{n=1}^N|d(n)|^2  .
\end{equation}
This inequality remains true without change when replacing the interval of integration by any other   of same length.   
\medskip\par
 Consider now   higher moments.
Let $q$ be some positive integer and denote 
$$E_q=\Big\{\underline{   k}=(k_1,\ldots, k_N); 
k_i\in \N : k_1+\ldots+ k_N= q\Big\}.$$ 
We assume that        
$ 
\l_1,
\ldots,
\l_N $ are     linearly independent  reals. The typical example is  $\l_j= \log p_j$, $j=1,\ldots,N$,
 where $ p_1, p_2,\ldots ,  p_N $ are different primes, see (\ref{liniprime}). 
\smallskip\par
Introduce a {\it coefficient   of linear spacing} of order
$q$  by putting
 $$ \xi_\l(N, q) =\inf_{\underline{   h}, \underline{   k}\in E_q\atop 
\underline{   h}\not= \underline{   k} } \big|(h_1 -k_1)\l_1+\ldots+   (h_N-k_N)\l_N \big|.  $$
     
   By assumption $\xi_\l(N, q) >0$ and
 $ \xi_\l(N, 1) =\inf_{1\le i,j\le N\atop 
i\not=j } \big| \l_i -\l_j \big| $.    If we expand the integrand, next integrate, we shall get similarly  
\begin{prop} For any interval $J$, denoting   
$|J|$ its length, \begin{equation}  {1\over |J|}\int_J \Big|  \sum_{n=1}^N   d(n)   e^{  it\l_n } \Big|^{2q}\dd t\le   
\big(\sum_{n=1}^N|d(n)|^{ 2 }\big)^{q} 
\Big( q! +{2 \min(N^{ q},\pi q!)\over |J|\xi}
\Big) .
\end{equation}
\end{prop} {\it Proof.} 
 Let
$J =[d,d+T]$. Put $P(t)= \sum_{n=1}^N   d(n)   e^{  it\l_n }$ and $\xi= \xi_\p(N, q)$. 
    Plainly 
 \begin{eqnarray}   \big|P(t ) \big|^{2q}\!\!   &=& \!\! \sum_{\underline{   k},\underline{   h}\in E_q}  {(q!)^2\over
k_1!h_1!\ldots k_N!h_N!} 
 \prod_{n=1}^Nd(n)^{ k_n}{\overline{d(n) } } ^{ h_n} e^{i t (k_n-h_n)\l_n }   \cr & =&
\sum_{\underline{   k} \in E_q } \Big({ q!  \over k_1! \ldots k_N! } 
\Big)^2\prod_{n=1}^N |d(n)|  ^{ 2k_n} +R (t)
\end{eqnarray} 
where 
\begin{equation}\label{3}R(t)  =  \sum_{\underline{   k},\underline{   h}\in E_q\atop 
 \underline{   k}\not=\underline{   h}} \Big({(q!)^2\over k_1!h_1!\ldots k_N!h_N!} 
\Big)\prod_{n=1}^N d(n) ^{ k_n} \overline{d(n) } ^{ h_n} e^{i t (k_n-h_n)\l_n  } .
\end{equation}  
 By integrating and using   linear independence
\begin{eqnarray}\label{ow}{1\over T}\int_J\big|P(t ) \big|^{2q}\dd t &=&\sum_{\underline{   k} \in E_q } \Big({ q! \over k_1! \ldots
k_N! } 
\Big)^2\prod_{n=1}^N     |d(n)| ^{ 2k_n}  
 \cr & & + \sum_{\underline{   k},\underline{   h}\in E_q\atop 
 \underline{   k}\not=\underline{   h}}  {(q!)^2\over k_1!h_1!\ldots k_N!h_N!} 
\prod_{n=1}^N d(n)  ^{ k_n} \overline{d(n) } ^{ h_n} 
\cr & &\times \Big[{  e^{i(d+ T)
\sum_{n=1}^N(k_n-h_n)\l_n }-   e^{i d 
\sum_{n=1}^N(k_n-h_n)\l_n }  \over iT(\sum_{n=1}^N(k_n-h_n)\l_n)}\Big] .
\end{eqnarray}   Put    
$${\bf c}_{\underline{   k}}= \prod_{n=1}^N{(d(n)  e^{i (d+T)   \l_n })^{ k_n}\over k_n!  }, \quad {\bf d}_{\underline{  
k}}= \prod_{n=1}^N{(d(n)  e^{i  d    \l_n })^{ k_n}\over k_n!  } ,\quad  {\bf   l}_{\underline{   k}}=\sum_{n=1}^N k_n \l_n .$$  
Then 
\begin{equation}\label{4} {1\over T}\int_J\big|P(t ) \big|^{2q}\dd t = (q!)^2
\sum_{\underline{   k} \in E_q } |{\bf d}_{\underline{   k}}|^2  
+{(q!)^2\over iT}\bigg\{\sum_{\underline{   k},\underline{   h}\in
E_q\atop 
 \underline{   k}\not=\underline{   h}}    { {\bf c}_{\underline{   k}}\overline{{\bf c}}_{\underline{   h}}\over  
{\bf   l}_{\underline{   k}}-{\bf   l}_{\underline{   h}}}    -\sum_{\underline{   k},\underline{   h}\in E_q\atop 
\underline{   k}\not=\underline{   h}}    { {\bf d}_{\underline{   k}}\overline{{\bf d}}_{\underline{   h}}\over   
{\bf   l}_{\underline{   k}}-{\bf   l}_{\underline{   h}}}    \bigg\} .
\end{equation}

We shall apply Hilbert's inequality under the following form:     let $\{x_{\underline{   k} } ,   y_{\underline{   k} }, \underline{  
k}\in E_q\}
$. Let also
$\{\l_{\underline{   k} },\underline{   k} \in E_q\} $  be distinct  real numbers such that $\min\{|\l_{ \underline{   k} }-\l_{
\underline{   h} }|,{ \underline{   k}}\not={
\underline{   h} }\}\ge \d$.
 Let   $\nu =\#\{E_q\}$ and consider  a bijection   
  $i: 
 \{ 1,\ldots,\nu\}\to E_q$. 
 By  using  Lemma \ref{HMV}  
\begin{eqnarray}\label{hima}\Big|\sum_{\underline{   k},\underline{   h}\in E_q\atop  \underline{   k}\not=\underline{  
h}}{x_{\underline{   k} }y_{\underline{   h} }\over
 \l_{ \underline{   k} }-\l_{
 \underline{   h} }}\Big| & =&\Big|\sum_{1\le u,v\le \nu\atop   u\not=v}{x_{i(u)
}y_{i(v) }\over
 \l_{ i(u) }-\l_{ i(v) }}\Big| \cr &\le& {\pi\over
\d}\Big(\sum_{1\le u \le \nu}  |x_{ i(u )}|^2\Big)^{1/2}\Big(\sum_{1\le  v\le \nu}  |y_{i(v) }|^2\Big)^{1/2}
\cr &= & {\pi\over
\d}\Big(\sum_{\underline{   k} \in E_q}  |x_{\underline{   k}  }|^2\Big)^{1/2}\Big(\sum_{  \underline{   h}\in E_q} 
|y_{\underline{  h} }|^2\Big)^{1/2}. 
\end{eqnarray}

Apply it to each of the two sums   in   parenthesis of the right-term  in (\ref{4}), we find  
\begin{equation}\label{hima} {(q!)^2\over  T}\bigg|\sum_{\underline{   k},\underline{   h}\in E_q\atop 
 \underline{   k}\not=\underline{   h}}    { {\bf c}_{\underline{   k}}\overline{{\bf c}}_{\underline{   h}}\over  
{\bf   l}_{\underline{   k}}-{\bf   l}_{\underline{   h}}}   -\sum_{\underline{   k},\underline{  
h}\in E_q\atop 
 \underline{   k}\not=\underline{   h}}    { {\bf d}_{\underline{   k}}\overline{{\bf d}}_{\underline{   h}}\over   
{\bf   l}_{\underline{   k}}-{\bf   l}_{\underline{   h}}}    \bigg|\le 
  {2\pi(q!)^2\over
  T\xi}\sum_{\underline{   k} \in E_q }      |{\bf d}_{\underline{   k}}|^2 \le  {2\pi q!\over
  T\xi}\,\Big[\sum_{n=1}^N |d(n) |^2\Big]^q,
\end{equation} 
since
\begin{eqnarray}\label{sa} (q!)^2\sum_{\underline{   k} \in E_q }      |{\bf d}_{\underline{   k}}|^2&=&\sum_{ k_1+\ldots+k_N=q } 
  \Big[ {  q!  \over k_1! \ldots k_N!
} \Big]^2   \prod_{n=1}^N|d(n) |^{ 2k_n}  \cr &\le &q!\,\sum_{ k_1+\ldots+k_N=q } 
    {  q!  \over k_1! \ldots k_N!
}    \prod_{n=1}^N|d(n) |^{ 2k_n}  = q!\,\Big[\sum_{n=1}^N |d(n) |^2\Big]^q .
\cr & &
\end{eqnarray}
 The   way to bound in (\ref{sa}), in turn, already appeared  in \cite{[S]}.

By substituting in (\ref{4}), we therefore  obtain
\begin{equation}{1\over T}\int_J\big|P(t ) \big|^{2q}\dd t \le     q!\Big(1+   {2\pi \over
  T\xi}
\Big) \,\Big[\sum_{n=1}^N |d(n) |^2\Big]^q.
\end{equation}
  Further,  from (\ref{ow}) we also get by using  Cauchy-Schwarz inequality 
  \begin{eqnarray}\label{ow1}{1\over T}\int_J\big|P(t ) \big|^{2q}\dd t &=&\sum_{\underline{   k} \in E_q } \Big({ q! \over
k_1!
\ldots k_N! } 
\Big)^2\prod_{n=1}^N     |d(n)| ^{ 2k_n}  
 \cr & & + \sum_{\underline{   k},\underline{   h}\in E_q\atop 
 \underline{   k}\not=\underline{   h}}  {(q!)^2\over k_1!h_1!\ldots k_N!h_N!} 
\prod_{n=1}^N d(n)  ^{ k_n} \overline{d(n) } ^{ h_n} 
\cr & &\times \Big[{  e^{i(d+ T)
\sum_{n=1}^N(k_n-h_n)\l_n }-   e^{i d 
\sum_{n=1}^N(k_n-h_n)\l_n }  \over iT(\sum_{n=1}^N(k_n-h_n)\l_n)}\Big]  \cr
   & \le  & q!\,\Big[\sum_{n=1}^N |d(n) |^2\Big]^q+{2\over T\xi}\Big(  \sum_{n=1}^N  |d(n) |\Big)^{2 q}
  \cr  &\le&   \Big( q! +{2 N^{ q}\over T\xi} \Big)\Big[\sum_{n=1}^N |d(n) |^2\Big]^q
  .
\end{eqnarray}
 Combining the two last estimates gives
\begin{equation}{1\over T}\int_J\big|P(t ) \big|^{2q}\dd t \le \Big( q! +{2 \min(N^{ q},\pi q!)\over T\xi}
\Big)\Big[\sum_{n=1}^N |d(n) |^2\Big]^q.
\end{equation}
 \cqfd 
\medskip\par
  Now we pass to high moments of  Dirichlet approximating polynomials
$${1\over T} \int_0^T \Big|  \sum_{n=1}^N  {1\over
n^{ {1\over 2}+it  } } \Big|^{2\nu}\dd t.$$ 
Apply Proposition \ref{MV} to   \begin{equation}\label{dir0}
 \Big(\sum_{n=1}^N  {1\over
n^{ {1\over 2}+it  } }\Big)^{\nu}:=\sum_{m=1}^{N^\nu}{b_m\over m^{ {1\over 2}+it  } }.
\end{equation}

Since $\d\ge 
 \min\big\{  \log(1+ {m-n\over n}) : 1\le n<m\le N^\nu\big\}\ge   {1\over 2 N^\nu}$, we
get  \begin{equation}\label{MV1}   \bigg|\int_0^T \Big|  \sum_{n=1}^N  {1\over
n^{ {1\over 2}+it  } } \Big|^{2\nu}\dd t -  T\sum_{m=1}^{N^\nu}{b_m^2 \over m  } \bigg|\le  CN^\nu \sum_{m=1}^{N^\nu}{b_m^2 \over m  }  .
\end{equation}
Recall that  $d_\nu(n)$ denotes the number of representations of the integer $n$ as a product of
$\nu$ factors. As    $b_m = \#\{ (n_1, \ldots, n_\nu); n_j\le N: m=n_1  \ldots n_\nu\} \le d_\nu(m)$, and  (\cite{Iv}, Section 9.5)
\begin{equation}\label{dm1}\sum_{m\le N} {d_{ \nu}^2(m)\over m}= (C_\nu +o(1))\log^{\nu^2} N,
\end{equation}
it follows that
 \begin{equation}  \label{asymp}  {1\over T} \int_0^T \Big|  \sum_{n=1}^N  {1\over
n^{ {1\over 2}+it  } } \Big|^{2\nu}\dd t \le C_\nu (1+{N^\nu\over T} ) \log^{\nu^2} N  .
\end{equation}  
 Hence if $T\ge N^\nu$
 \begin{equation}  \label{mvnu}     {1\over T} \int_0^T \Big|  \sum_{n=1}^N  {1\over
n^{ {1\over 2}+it  } } \Big|^{2\nu}\dd t \le C_\nu   \log^{\nu^2} N  .
\end{equation}
The latter estimate is in fact two-sided, see Corollary \ref{10}. It can also be reformulated as 
 \begin{equation}  \label{mvnu0}  c_\nu  \sum_{m=1}^{N^\nu} {b_m^2\over m} \le    {1\over T} \int_0^T \Big|  \sum_{n=1}^N  {1\over
n^{ {1\over 2}+it  } } \Big|^{2\nu}\dd t \le C_\nu  \sum_{m=1}^{N^\nu} {b_m^2\over m}  .
\end{equation}
Now,   
  by  using    approximation formula (\ref{approx}),   it follows that
the   reformulation (\ref{equivdemi}) of the LH is also equivalent to  \begin{equation} \label{approx4}{1\over  N}\int_{0}^{N} \Big| 
\sum_{n=1 }^N  {1\over  n^{ {1\over 2}+it}} -{N^{{1\over 2}-it}\over {1\over 2}-it}   \Big|^{ 2\nu} \dd t  = \mathcal O_\e(N^\e) \quad\qq
\nu=1,2,\ldots.
\end{equation}
The critical range of values of $T$  in (\ref{mvnu}) is thus $T\sim N $. But   it is a simple matter to observe   that in this
case, estimate (\ref{mvnu}) can no  longer 
 be true, unless the LH is false.  
 Indeed by the Minkowski inequality, if (\ref{mvnu}) and (\ref{approx4}) were simultaneously true, we would have
 \begin{equation}  {1\over  N}\int_{0}^{N} \Big| {N^{{1\over 2}-it}\over {1\over 2}-it} \Big|^{ 2\nu} \dd t  ={N^{ \nu-1}\over 2 }
\int_{0}^{N}
  {  \dd t\over ({1\over 4}+t^2)^\nu} \sim  C N^{ \nu-1}   = \mathcal O_\e(N^\e),
\end{equation}
which is absurd as soon as $\nu>1$. 
 It also follows from these observations that the order of 
$${1\over T} \int_0^T \Big|  \sum_{n=1}^N  {1\over
n^{ {1\over 2}+it  } } \Big|^{2\nu}\dd t$$
  is necessarily much bigger for $T\le N^\nu$ than for $T\ge N^\nu$.
\bigskip\par   
  Concerning upper bounds,  there is a useful majorization argument (\cite{[Mon]}, p.131) which can be stated for arbitrary even powers. 
 \begin{prop} Let $q$ be
any positive integer. Let     
$c_1,\ldots, c_N$ be   complex numbers and nonnegative reals $a_1,\ldots, a_N$  such that $|c_n|\le a_n$, 
$n=1,\ldots, N$. Then for any reals $ \p_1,
\ldots,
\p_N  $ and any reals  $ T,T_0$ with
$T>0$ 
\begin{equation} \label{ma}\int_{|t- T_0|\le T}\Big| \sum_{n=1}^N c_ne^{it\p_n}\Big|^{2q}dt \le 3 \int_{|t |\le T}   \Big|
\sum_{n=1}^N a_ne^{it\p_n}\Big|^{2q}dt. 
\end{equation}
\end{prop}

 Consider the kernel 
 $$ K_T(t)=K_T(|t|)=\big( 1-|t|/T)\chi_{\{|t|\le T\}}   . $$
 The proof of (\ref{ma}) is based on the following   properties of $K_T$: for any reals $t,H$
\begin{eqnarray*}
a)\quad&&\chi_{\{|t-H|\le T\}} \le K_T(t-H)+K_T(t-H+T)+ K_T(t-H-T)\cr  
b)\quad &&\widehat{ K}_T(u)  =T^{-1}\big({\sin  Tu\over   u  }\big)^2  \qq\hbox{for any real $u$}  .
\end{eqnarray*}

From (\ref{ma})   one can derive the following lower bound \cite{W6}. 
\begin{thm}   For any positive integer $q$,
there exists a constant $c_q$, such that for any reals $ \p_1,
\ldots,
\p_N  $,   any non-negative reals   $a_1,\ldots, a_N$, and any $T>0$,    
\begin{equation}  \label{t3} c_{q}\Big(\sum_{n=1}^N  a_n^2\Big)^q  \le {1 \over 2T} \int_{|t |\le
T}   \Big|
\sum_{n=1}^N a_n e^{it\p_n}\Big|^{2q}dt. 
\end{equation}
\end{thm}   The $L^1$-case is related to   Ingham's inequality. Recall the sharp  form   due to
Mordell \cite{[Mor]}: let
$0<\p_1<\ldots<\p_N$ and  let $\gamma$ be such that
$\displaystyle{\min_{1<n\le N}} \p_n-\p_{n-1}\ge \gamma>0$. Then  
\begin{equation}\sup_{n=1}^N|  a_n| \le {K\over T}\int_{-T}^T\Big|
\sum_{n=1}^N a_n e^{it\p_n}\Big| dt\qq {\rm with}\ T={\pi\over \gamma}, \label{4.1}\end{equation} 
where $K\le 1$.  
 Further  with no restriction, one
always have the very familiar inequality    in the theory of uniformly almost periodic functions:
\begin{equation}\sup_{n=1}^N |a_n|\le \limsup_{T\to \infty}{1\over 2T}\int_{-T}^T\Big|
\sum_{n=1}^N a_n e^{it\p_n}\Big| dt \le \sup_{t\in \R}\Big|
\sum_{n=1}^N a_n e^{it\p_n}\Big|, \label{4.2}\end{equation}
   Inequality  (\ref{t3})  is obtained    by choosing    $c_n=\e_na_n $ in (\ref{ma}), where
$\ue=\{\e_n, n\ge 1\}$ is a Rademacher sequence, next taking expectation   and using   Khintchine-Kahane inequalities for Rademacher sums.
   We shall deduce from it  the following lower bound.
\begin{cor}\label{10}
 {\it  For every $N$, $T$ and $\nu$} $$c_\nu\log^{\nu^2}  N\le
{1
\over 2T}
\int_{| t |\le T}   \Big|
\sum_{n=1}^N {1\over n^{{1\over 2}+i t}} \Big|^{2\nu}d t. $$
\end{cor}

 Indeed, apply (\ref{t3})  with   $q=2$ to  the sum
$$
 \Big(\sum_{n=1}^N  {1\over
n^{ {1\over 2}+it  } }\Big)^{\nu}:=\sum_{m=1}^{N^\nu}{b_m\over m^{ {1\over 2}+it  } }.
$$
  Then for all $N$ and $T$
 $$c_\nu\ \sum_{m=1}^{N^\nu}{b^2_m\over m  } \le
{1
\over 2T}
\int_{| t |\le T}   \Big|
\sum_{n=1}^N {1\over n^{{1\over 2}+i t}} \Big|^{2\nu}d t. $$
Notice that if $m\le N$, $b_m=d_{ \nu}(m)$ and we know that
$$\sum_{m\le x} {d_{ \nu}^2(m)\over m}= (C_\nu +o(1))\log^{\nu^2} x. $$
See  
\cite{Iv} section 9.5. Thus 
$$\sum_{m=1}^{N^\nu}{b^2_m\over m  }\ge \sum_{m=1}^{N }{b^2_m\over m  } \ge c_\nu \log^{\nu^2} N$$
Henceforth $$c_\nu\log^{\nu^2}  N\le
{1
\over 2T}
\int_{| t |\le T}   \Big|
\sum_{n=1}^N {1\over n^{{1\over 2}+i t}} \Big|^{2\nu}d t. $$
\cqfd 
 
    Put 
 $$ M_\nu (T)=   \int_{0}^{T}   \big|
\zeta({1\over 2}+it ) \big|^{2\nu}dt$$
  The corresponding inequality for the Riemann Zeta-function is Ramachandra's  well-known lower bound and we
recall   (\cite{Ti} p.180, see also
\cite{Iv} section 9.5 and \cite{Ra})  that
\begin{equation} \label{ram}   c_{\nu}T(\log   T)^{\nu^2} \le M_\nu (T). 
\end{equation}
    
Assuming RH,  Soundararajan recently proved in \cite{So} that for every positive real number $\nu$, and every$\e>0$, we have
\begin{equation} \label{soun}       M_\nu (T)\le C_{\nu, \e}\ T(\log   T)^{\nu^2+\e}. 
\end{equation} We conclude this section by
mentionning and briefly discussing some    related   results for the Riemann
$\zeta$-function.   

\begin{rem} \rm {\gsec (Mean value results for the   $\zeta$-function)}     For the critical value $\s=1/2$, the most achieved results are
\begin{equation}\label{mv2zeta}\int_0^T |\zeta({1\over 2} +it)|^2 dt =T\log\Big({T\over 2\pi} \Big) + (2\gamma -1)T +E(T),
\end{equation}
where $\gamma$ is Euler's constant and the error term $E(T)$ satisfies $E(T)\ll_\e T^{1/3 +\e}$, see \cite{Ti} p.176. And
\begin{equation} \label{mv4zeta}\int_0^T |\zeta({1\over 2} +it)|^4 dt ={T\log^4 T  \over 2\pi^2} +  \mathcal O(T\log^3 T),
\end{equation}
see \cite{Ti} p.148. The approximate equation (\ref{approx}) already suffices to show
\begin{equation} \int_0^T |\zeta({1\over 2} +it)|^2 dt =\mathcal O(T\log T) .
\end{equation}
The very formulation of (\ref{approx}) yields for the fourth moment that it is equivalent to work with 
$$\Big|\sum_{n\le x} {1\over n^s} -{x^{1-s}\over 1-s}\Big|^4 , \qq x\sim T,$$
instead of $|\zeta({1\over 2} +it)|^4$, when $0\le t\le T$. However there is apparently no  known proof   of  $\int_0^T |\zeta({1\over 2}
+it)|^4 dt = 
\mathcal O(T\log^4 T)$ based on (\ref{approx}), which is a bit frustrating. 
 In place, one has to use the following more elaborated 
approximate equation 
\begin{equation} \label{approx1/2} \zeta(s)= \sum_{n\le x} {1\over n^s} +\chi(s)  \sum_{n\le y} {1\over n^{1-s}} +\mathcal O( x^{-\s}\log |t|)+\mathcal
O(|t|^{{1\over 2}-\s} y^{ \s-1}),
\end{equation}
in which $h$ is a positive constant, $0<\s<1$, $2\pi xy=t$, $x>h>0$, $y>h>0$ and 
$$\chi (s)  = {\zeta(s)\over \zeta(1-s)}=2^{s-1} \pi^2 \sec \Big({s\pi\over 2\Gamma (s)} \Big). $$
This function verifies in any fixed strip $\a\le \s\le \b$,  
 $|\chi(s)|\sim \big({t/ 2\pi}\big)$,  
as $t\to \infty$.
\medskip\par
  The knewledge concerning     moments 
$$\int_0^T |\zeta({1\over 2} +it)|^k \dd t 
$$ beyond  $k=4$ is,    to say the least, very sparse.
 For the case $k=12$, we may quote the beautiful   result due to Heath-Brown
$\int_0^T |\zeta({1\over 2} +it)|^{12} dt \ll T^2 \log^{17} T$. See  \cite{Ti} p. 79, 95  and 178 for the aforementionned facts. See also
\cite{Iv} Section 8.3.
 
There are also alternative  mean-value theorems involving integrals of the form
$$ J(\d)= \int_0^\infty |\zeta({1\over 2} +it)|^{2k} e^{-\d t} \dd t, \qq \d\to 0. $$
The behavior of these integrals is similar to the one of 
$$ I(T)= \int_0^T |\zeta({1\over 2} +it)|^{2k}   \dd t, \qq T\to \infty, $$
and this follows notably from integral versions of a well-known Tauberian result of Hardy and Littlewood. More precisely, if $ f\ge 0$, then (\cite{Ti}
Ch. VII)
\begin{equation}  \int_0^\infty f(t) e^{-\d t} \dd t\ \buildrel{\d \to 0}\over{ {}_{\widetilde{\quad} }}{1\over \d} \qq \Rightarrow \qq \int_0^T f(t) 
\dd t
\ \buildrel{T \to \infty}\over{ {}_{\widetilde{\quad} }}T.\end{equation}
When $1/2<\s <1$, we have   
\begin{equation}\label{mv2}\int_1^T|\zeta( \s +it)|^2 \dd t =T \sum_{n=1}^\infty {1\over n^{ 2\s}} + \mathcal
O  \big( T^{\,  2 -2\s  }\big).\end{equation}
\begin{equation}\label{mv4}\int_1^T|\zeta( \s +it)|^4 \dd t =T \sum_{n=1}^\infty {d_2^2(n)\over n^{ 2\s}} + \mathcal
O_\e \big( T^{\, 3/2 -\s +\e}\big).
\end{equation}
 
There are also results for $k=1/2$. For   $k>2$ integer, it is known (\cite{Ti},  p.125)
that
\begin{equation}\lim_{T\to
\infty}{1\over T}\int_1^T
\big|\zeta(\s + it)\big|^{2k} {\rm d}t=
\sum_{n=1}^\infty {d_k^2(n)\over n^{2\s}}, \qq (\s>1-1/k).\end{equation}
We refer to \cite{Iv}, Chapter 8 and notably Theorem 8.5 for improvments of this,  up to power $12$, under weaker   conditions on $\s $.
 \end{rem}
\begin{rem}\rm   {\gsec (Square function of the Riemann-zeta function)}   Let
$\theta =\{T_j, j\ge 1\}$ be such that
$T_j\uparrow\infty$.  Given any fixed positive integer
$k$,  we  define for $1/2<\s<1$ the $\zeta$-square function $\mathcal S_\theta(\zeta, \s)$ associated to $\theta$ as follows
$$ \mathcal
S_\theta(\zeta, \s):=\bigg(\sum_{j=1}^\infty  \Big|{1\over T_{j+1}}\int_0^{T_{j+1}}|\zeta(\s +it)|^{2k}dt -{1\over
T_{j}}\int_0^{T_{j}}|\zeta(\s +it)|^{2k}\Big|^2\bigg)^{1/2}.$$   The finiteness in (\ref{sqfif})  of the square function linked to the
Fourier inversion formula   and the analogy   described at the beginning of Section 2 between Lindel\"of Hypothesis and Fourier inversion
formula suggest  to
  investigate properties of  the
$\zeta$-square function $\mathcal S_\theta(\zeta, \s)$. 
\smallskip\par\noi 
  {\gsec Problem.} For which sequences $\theta$ is $\mathcal S_\theta(\zeta, \s)$ finite for all $1/2<\s<1$? When is the same   also
true  independently from the value of
$k\ge 1$?
\vskip 2pt
In the case $k=1$, $k=2$, it follows trivially from (\ref{mv2}), (\ref{mv4}) that any geometrically increasing sequence $\theta$,  
$T_{j+1}/ T_j
\ge M>1$, is suitable. One may wonder whether this condition is also necessary for the finiteness of $\mathcal
S_\theta(\zeta,
\s)$  for every
$1/2<\s<1$. The case $k=1/2$ is also of interest.


 \medskip\par  At this regard, it is worth noticing that   for any integers $n\ge m\ge 1$, we have
\begin{eqnarray}\label{f3}  & &\!\!\!\!\! \!\!\!\!\!\!\!\!\!\!\!\!\!\!\!\!\!\!\!\!\Big|\sum_{k=1}^n   {1\over k^{ \frac12+it}}
-{n^{\frac12-it}\over
 \frac12-it}\Big|^2-\Big|\sum_{k=1}^m   {1\over k^{ \frac12+it}} -{m^{\frac12-it}\over
 \frac12-it}\Big|^2
\cr &=& \ \ \!\!\!\!\!\!\!\!\sum_{\ell=m+ 1}^n{1\over \ell }    +2 \sum_{m+ 1\le k\le n\atop 1\le \ell<k}
  \Re\bigg\{{1\over \sqrt{k\ell}}e^{it \log {k\over \ell}}- \int_{\log {k-1\over \ell}}^{\log {k\over \ell}}
  e^{({1\over 2} +it)x} \dd x \bigg\}.
\end{eqnarray} This  follows   from a more homogeneous
reformulation
 of  the approximating term in (\ref{approx}),  namely
 \begin{equation}\label{f1}\Big|\sum_{k=1}^n   {1\over k^{ \frac12+it}} -{n^{\frac12-it}\over
 \frac12-it}\Big|^2=\sum_{\ell= 1}^n{1\over \ell }    +2 \sum_{1\le \ell<k\le n}
  \Re\bigg\{{1\over \sqrt{k\ell}}e^{it \log {k\over \ell}}- \int_{\log {k-1\over \ell}}^{\log {k\over \ell}}
  e^{({1\over 2} +it)x} \dd x \bigg\}.
\end{equation}  
  Indeed,
\begin{eqnarray*}   &  &  \sum_{\ell= 1}^n{1\over \ell }    +2 \sum_{1\le \ell<k\le n}
  \Re\bigg\{{1\over \sqrt{k\ell}}e^{it \log {k\over \ell}}- \int_{\log {k-1\over \ell}}^{\log {k\over \ell}}
  e^{({1\over 2} +it)x} \dd x \bigg\}\cr &=& \sum_{\ell= 1}^n{1\over \ell }    +2 \sum_{1\le \ell<k\le n}
  \Re\bigg\{{1\over \sqrt{k\ell}}e^{it \log {k\over \ell}} \bigg\} -2 \sum_{1\le \ell<k\le n}
  \Re\bigg\{  \int_{\log {k-1\over \ell}}^{\log {k\over \ell}}
  e^{({1\over 2} +it)x} \dd x \bigg\}
\cr &=&  \Big|\sum_{k=1}^n   {1\over k^{ \frac12+it}} \Big|^2 -2 \sum_{ \ell=1}^{n-1}
  \Re\bigg\{  \int_{0}^{\log {n\over \ell}}
  e^{({1\over 2} +it)x} \dd x \bigg\}
\cr &=&  \Big|\sum_{k=1}^n   {1\over k^{ \frac12+it}} \Big|^2 -2 \sum_{ \ell=1}^{n-1}
  \Re\bigg\{  
  {1\over {1\over 2} +it}\Big[\Big({n\over \ell}\Big)^{{1\over 2} +it}  -1\Big]\bigg\}  
\cr &=&  \Big|\sum_{k=1}^n   {1\over k^{ \frac12+it}} \Big|^2 -  2 \Re\bigg\{ {  n^{{1\over 2} +it}\over ({1\over 2} +it) } 
\sum_{ \ell=1}^{n-1}
   {1\over \ell^{{1\over 2} +it}}   \bigg\}+   
    {n  \over {1\over 4} + t^2} . 
\end{eqnarray*}

 If $1/2< \s<1$,  there is  a similar  formula (\cite{W5},
Corollary 5)
   \begin{eqnarray}  \label{f2}\!\!\!\!\!\!\!\!\!\!\ \Big|\sum_{k=1}^n   {1\over k^{ \s+it}} -{n^{\s-it}\over
 \s-it}\Big|^2    & =&  \sum_{\ell =  1}^n{1\over \ell^{2\s}}+2 \sum_{1\le \ell<k\le n}
  \Re\bigg\{{e^{it \log {k\over \ell}}\over (k\ell)^\s}  \cr&&\!\!\!\!\!\!\!\!\!\!\!- \ell^{1-2\s}\int_{\log {k-1\over \ell}}^{\log
{k\over
\ell}}
  e^{(1-\s +it)x} \dd x \bigg\}-    { \Psi_\s \over  (1-\s)^{ 2}+  t^2
 }, 
\end{eqnarray} 
where   $\ \displaystyle{\Psi_\s=  \s + (1-2\s)  2\s  
        \sum_{k=1}^{\infty} \int_0^1
        \frac {t-t^2}{2} (k+t)^{  -2\s-1} dt   + \mathcal O(n^{1-2\s}).
 } $ 
\end{rem}
 
  \section{Supremum of Dirichlet polynomials}   
     We begin with some general considerations.
 Let $d:\N\to \R$.
 The   supremum  of the Dirichlet polynomials
 $ P(s)= \sum_{n=1}^N  d(n) n^{-s}
 $
over lines $\{s=\s+it, \ t\in \R\}$ is naturally related to that 
of corresponding Dirichlet series, via the abscissa of uniform convergence
$$ \s_u=\inf\Big\{\s : \sum_{n=1}^\infty  d(n) n^{-\s - it}\
\hbox{converges uniformly over $t\in\R$} \Big\},
$$
through the relation
\begin{equation}\s_u= \limsup_{N\to \infty} {{\log  \,  \displaystyle
\sup_{t\in \R} \big|\sum_{n=1}^N  d(n) n^{-it}\big|} \over \log N}\  .
\end{equation}     We   refer to   \cite{B},\cite{H} or 
\cite{HR} for   background
and related results. This naturally   justifies the
 investigation of the supremum of Dirichlet polynomials.

  \medskip \par
It is natural to first compare the behavior of the suprema of Dirichlet polynomials with the one of trigonometric polynomials, which we
shall do by investigating Rudin-Shapiro polynomials.  We refer to \cite{[Mon]} Chapter 7 where a comparative study is presented.
\medskip \par\noi   {\gsec Rudin-Shapiro polynomials.}
 Recall   the classical
setting. For any trigonometric polynomial we have
\begin{equation}
{\sum_{n=1}^{N }
|d(n)| \over \sqrt{N}} \le
\sup_{t\in \R} | \sum_{n=1}^{N } d(n) e^{int} |
\le 
\sum_{n=1}^{N } |d(n)|.
\end{equation}
The arguments for getting the
lower bound are the inequality between the sup-norm and
$L_2$-norm, the
orthogonality of $(e^{int})_n$ and H\"older inequality. 

Rudin and Shapiro
constructed a fairly simple sequence $d(n)\in\{-1,+1\}$ such that the
right
order of the lower bound is attained:
\begin{equation}\sup_{t\in \R} |
\sum_{n=1}^{N } d(n) e^{int} |
\le (2+\sqrt{2}) \sqrt{N+1}
\sim
(2+\sqrt{2}) \ {\sum_{n=1}^{N } |d(n)| \over \sqrt{N}}.
\end{equation}
Consider
now the Dirichlet polynomials instead of the trigonometric ones. It is
known from
\cite{KQ} and \cite{Q3} that 
\begin{thm} \label{supKQ}For any $(d(n))$
\begin{equation}
\sup_{t\in \R} |
\sum_{n=1}^{N } d(n) n^{it} | \ge 
\alpha_1 {\sum_{n=1}^{N } |d(n)| \over
\sqrt{N}} 
\exp\{\beta_1 \sqrt{\log N\log\log N}\}.
\end{equation}
and for some
$(d(n))$
\begin{equation} 
\sup_{t\in \R} | \sum_{n=1}^{N } d(n) n^{it} | \le 
\alpha_2
{\sum_{n=1}^{N } |d(n)| \over \sqrt{N}} 
\exp\{\beta_2 \sqrt{\log N\log\log
N}\},
\end{equation}
with some universal constants
$\alpha_{1 },\alpha_{ 2},\beta_{1 }, \beta_{ 2}$.
\end{thm}
  A finer result    with explicit constants was  recently   obtained by de la Bret\`eche in
\cite{Br}.
 Therefore  the lower bound for
Dirichlet polynomials is necessarily worse  than in the classical case. Notice also
that the construction
  in \cite{Q3} is a probabilistic one; no
explicit example 
of Rudin-Shapiro type is known for Dirichlet
polynomials.
\medskip\par
\smallskip\par There is a basic reduction step in the study of the suprema. Introduce a useful notion. A set of numbers
$\p_1,\p_2,\ldots ,
\p_k
$ is linearly independent if no linear relation
 $ a_1\p_1+a_2\p_2+\ldots +a_k\p_k=0$, with integral coefficients, not all zero, holds between them.
 For a proof of the classical result below,   we refer to \cite{HW}. \smallskip\par\noi
{\gsec Kronecker's theorem.} If
$\p_1,\p_2,\ldots ,
\p_k,1$ are  linearly independent,
$\t_1,\t_2,\ldots ,
\t_k
$ are arbitrary, and $N$, $\e$ are positive, then there are  integers $n>N$, $n_1,n_2, \ldots, n_k $ such that 
$$ \max_{1\le m\le k} |n\p_m-n_m -\t_m|<\e.$$
 Consequently, the set of points 
 $\{n\p_1 \}, \{n\p_1\} , \ldots, \{n\p_k \}  $
is dense in $\T^k$. 
\smallskip\par
Let $ p_1, p_2,\ldots ,  p_k $ be different primes. By the fundamental theorem of arithmetic 
 \begin{equation}\label{liniprime}\log p_1,\log p_2, \ldots, \log p_k\qq \hbox{are linearly independent.}
\end{equation}  This will    enable to replace the Dirichlet polynomial by some relevant {\it trigonometric} polynomial.
Introduce the necessary notation. Let $2=p_1<p_2<\ldots$ be the sequence of consecutive primes. If $  n=\prod_{j=1}^\tau p_j^{a_j(n)}$, we
 write
$\underline{a}(n)=\big\{a_j(n), 1\le j\le \tau\big\}$. According to the standard notation we also denote   $\O(n)=a_1(n)+\ldots+a_\tau(n)$ and by $\pi(N)$ 
the number of prime numbers   less or equal to
$N$. Let us fix $N$.
We put in what follows $\m =\pi(N)$  and define  for $\underline{z}= (z_1,\ldots, z_\m) \in
\T^\m$,
$$ Q(\underline{z})= \sum_{n=2}^N
 {d}(n) n^{-\s}e^{2i\pi\langle \underline{a}(n),\underline{z}\rangle} ,
$$
  H. Bohr's observation 
states that
 \begin{equation}\label{bohr}\sup_{t\in \R} \big|\sum_{n=2}^N  d(n) n^{-(\s +it)}
 \big| =\sup_{\underline{z} \in \T^\m}
\big|Q(\underline{z})\big|\ .
\end{equation}
\begin{rem} \rm  Naturally  no similar  reduction   occurs  when     considering    the supremum over a given bounded
interval $I$. 
  However, when the length of $I$ is of
 exponential size with respect to the degree of $P$, precisely when 
 $$|I|\ge  e^{    (1+\e) \o  N(\log  N\o) 
  \log N  },$$  
  the related supremum   becomes
comparable, for $\o$ large,  to the one taken on the real line, with an error term of order $\mathcal O(\o^{-1}) $. This is in turn a
rather general phenomenon due to   existence of "localized" versions of Kronecker's theorem;    and in the present case to
   Tur\'an's estimate (see  
\cite{W3} for   a   slighly improved form of it using a probabilistic approach, and
      references therein).
When the length is of sub-exponential order, the study however still  belong to the field of
application of the general theory of regularity of stochastic processes. 
\end{rem}  

  Before going further notice, as an immediate consequence of Kronecker's Theorem, that if $\p_1,\p_2,\ldots ,
\p_N
$ are linearly independent then
\begin{equation}    \sup_{t\in \R} \big|\sum_{n=1}^N  d(n) e^{- it\p_n}
 \big| =\sum_{n=1}^N  |d(n)|.
 \end{equation}
 
Let us first consider lower bounds. Subsets  $A \subseteq\{1, \ldots, N\}$ such that 
$$\forall \{\d_n, n\in A\}\in  \big\{0,\frac12 \big\}^A,\quad \exists \uz \in \T^\tau :     \sum_{j=1}^\tau a_j(n)z_j= \d_n \   {\rm mod}(1),\quad \forall
n\in A 
$$ 
   are of particular interest, since   $e^{2 i \pi\langle \ua(n),\uz\rangle}=1$ or $-1$ according to $\d_n=0$ or $1/2$.
  As   $ \langle \ua(p_j),\uz\rangle=   z_j$, by choosing $\uz$ so that $z_j=0$ or $ 1/2$,   we deduce with (\ref{bohr})  
 \medskip\par\noi {\gsec Bohr's
 lower bound}  (\cite{Bo}) 
$$ \E\sup_{t\in \R} \big|\sum_{n=2}^N  d(n) n^{- it }
 \big|  \ge  \sum_{p\le N  } |d(p)|  p^{-\s}.$$
 This was generalized in \cite{Q2}  by Queff\'elec who proved 
\begin{prop}  For any integer $m\ge 1$
$$ C_m\,  \E\sup_{t\in \R} \big|\sum_{n=2}^N  d(n) n^{- it }
 \big|  \ge  \Big(\sum_{   n\le N\atop \O(n)=m}|d(n)|^{{2m \over m+1}}\Big)^{{m+1\over 2m}},$$
where $C_m=\big({2\over \sqrt \pi}\big)^{m-1}{m^{{m\over 2}}(m+1)^{{m+1\over 2}}\over 2^{m} (m!)^{{2\over m+1}}} $, ($C_1=1$). Further $C_m\le
m^{{m\over 2}}$.\end{prop}
These  estimates are crucial (\cite{KQ}, see section 4) in the proof of Theorem \ref{supKQ}. 
\medskip \par

\medskip\par\noi{\gsec Local suprema of Dirichlet polynomials.}   Let        
$ 
\p_1,
\ldots,
\p_N $ be     linearly independent  reals.   In \cite{W7}, the local suprema of the Dirichlet polynomials
   $P(t)= \sum_{n=1}^N c_ne^{it\p_n}$ is investigated. 
  Let $q$ be some positive integer.
    Then (\cite{W7}, Theorem 4), 
\smallskip\par\noi {\it There exists a constant
$C_q$ depending on $q$ only, such that for any intervals $J, L$  
 \begin{eqnarray} \Big({1\over |J|}\int_{J } \big |\sup_{t \in L}| P( \t+t)|d\t 
\Big)^{1/2q} & \le& C_q\, \mathcal B \max\Big\{1, |L|\tilde{ \p}_N \Big\}^{1/ 2q}
\bigg\{\Big[\sum_{n=1}^N |c_n|^2 
\Big]^{1/ 2} +  \cr &{}&      \min\Big(|L|,  {1\over \tilde{ \p}_N} \Big)   \Big[\sum_{n=1}^N |c_n|^2 \p_n^2
\Big]^{1/ 2}\bigg\} , 
\end{eqnarray}
 where 
  $\mathcal B =\big[q!  \big( 1 +{2  \pi  \over |J|\xi_\p(N,q) } \big)\big]^{1/2q}$, $\tilde{
\p}_N=
\sup_{n\le N} |\p_n|$}.

 This result is used in the same paper to investigate  by means of Tur\'an's result (\ref{zfr}), zerofree regions of the Riemann-zeta
function.
  \section{Random Dirichlet polynomials}   Studies   for
  random  Dirichlet polynomials and    random 
Dirichlet series were   developed in   \cite{Ha1} and    \cite{Q1},\cite{Q2},\cite{Q3} notably, see also  \cite{LW1},\cite{LW2} and
references therein.    Such   investigations concerning   random  Dirichlet series and random 
power  series  go back to earlier
works of Hartman \cite{Har}, Clarke \cite{C} and  Dvoretzky-Erd\" os   \cite{DE1},\cite{DE2}.
 
 Let us first quote some   general results. For instance let $ \underline \xi=\{\xi, \xi_n, n\ge 1\}$ be a 
sequence of i.i.d. random variables and let   $\sigma_c$ and 
$\sigma_a$ be, respectively, the
almost sure abscissa of convergence and of absolute convergence 
of the Dirichlet series $\sum_{n=1}^\infty \xi_nn^{-s}$. 
If $\xi\neq 0$ holds
with positive probability, let 
$k_\xi :=\sup\{\gamma : \E |\xi|^\gamma <\infty\}$. 
The connection between the abscissas $\sigma_c$ and $\sigma_a$ 
and integrability of $\xi$ has been clarified by Clarke 
in  
 \cite{C}. 
\begin{prop}We have the implications:
\begin{equation}  \matrix { k_\xi=0&\Rightarrow &\sigma_a=\sigma_c=\infty \cr
        0<k_\xi\le  1 & \Rightarrow &\sigma_a=\sigma_c=1/k_\xi  \cr      
(k_\xi>1 \ and \ \E\xi\neq 0)&\Rightarrow &\sigma_a=\sigma_c=1     \cr   
   (k_\xi>1 \ and \ \E\xi=0)&\Rightarrow &\sigma_a=1     \ and \   \sigma_c=\max(1/k_\xi,{  {1/2}})  .   }   
                                                   \end{equation} 
\end{prop} \medskip \par Now let here and throughout the remainding part of the paper
$\e= \{\e_i, i\ge 1\}$  be 
a sequence of independent Rademacher random variables
($\P\{ \e_i =\pm 1\} =1/2$) with basic
probability space $(\O, \A, \P)$.  The following 
result is due to  Bayart, Konyagin and Qu\'effelec \cite{BKQ}. 

\begin{thm} Let $\{d(n), n\ge 1\}$ be a sequence of complex numbers. If 
$$\limsup_{N\to \infty} {1\over \log \log N}
\sum_{n=0}^N |d(n)|^2=\g>0,$$
 then
 for almost all $\o$   the series $
\sum_{n=0}^\infty 
\e_n(\o) d(n) n^{i t} $ diverges for each $t\in \R$.
\end{thm}
  The result is nearly optimal: if $0<\d_n\to 0$, there 
exists a sequence $\{d(n), n\ge 1\}$ such that $\limsup_{N\to \infty} 
{1\over \d_N \log \log
N}\sum_{n=0}^N |d(n)|^2 >0$, but for each $\o$, the series $
\sum_{n=0}^\infty 
\e_n(\o) d(n) n^{i t}  $ converges for at least on $t\in \R$.
 \smallskip\par
 In relation with the above, we may quote Hedenmalm and 
Saksman's  extension \cite{HS}  of Carleson's result:
\begin{thm} Under the assumption $ {\sum_{n=0}^\infty |d(n)|^2<\infty}$
 the Dirichlet series 
  $$
\sum_{n=0}^\infty 
\e_n  d(n) n^{-1/2+i t}  $$
converges for almost all $t$. 
\end{thm}
 A simple and elegant proof is given in  Konyagin and Qu\'effelec \cite{KQ} p.158/159.

  \section{Supremum of random Dirichlet polynomials}  Now consider the random Dirichlet polynomials
 
\begin{equation} \label{e12}
  \DD(s) =\sum_{n=1}^N \e_n d(n) n^{-s}, \qq s=\s +it,
\end{equation}   
and examine their    supremum properties. When ${d}(n)\equiv 1$, there are   optimal results. If $\s=0$, then for some absolute constant
$C$, and all integers
$N\ge 2$
\begin{equation}  C ^{-1} {N  \over   \log  N }\le
\E\, \sup_{t \in \R}|\sum_{n=2}^N \e_n  n^{  - it}|
\le  C {N\over \log N}\ .
\end{equation}
This has been proved by Hal\'asz and  was later extended  by Queff\'elec  
to the range of values
$0\le \s<1/2$. Queff\'elec provided a probabilistic proof of the
original one, using Bernstein's inequality for polynomials.
\begin{thm} There exists a constant $C_\s$ depending on $\s$
only, such that  for   all integers $N\ge 2$  we have,
\begin{equation}\label{q1} C_\s ^{-1} {N^{1-\s}  \over   \log  N }\le \E\,
\sup_{t \in \R}|\sum_{n=2}^N \e_n  n^{-\s  - it}|  \le    C_\s
{N^{1-\s} \over \log  N }  .
\end{equation}
\end{thm}
Extensions of (\ref{q1})  were obtained in the recent works \cite{LW1},\cite{LW2}. The approach used does not appeal to  Bernstein's
inequality, and is completely based on stochastic process method, notably the metric entropy method. Further a new lower bound is obtained,
which is of a completely different nature than Bohr's deterministic lower bound used in Queff\'elec's proof.  
 For random Dirichlet polynomials defined in (\ref{e12}), a new approach is developed in \cite{LW1}. Define
 $$ {\cal L}_j=\Big\{n=p_j \, \tilde n\ : \
  \tilde n\le {N\over p_{j}}\ \hbox{and}\
P^+(\tilde n)\le p_{\mu/2}\Big\}, \qq  \qq j\in(\mu/2,\mu]. $$
   \begin{thm} \label{lb} {\gsec (Lower bound)}
$$\E \sup_{t\in \R} |\mathcal D (\s+it) | \ge c\ \sum_{ \mu/2<j\le \mu }
\Big(\sum_{n\in  {\cal L}_j}
d(n)^2\,  n^{-2\s} \Big)^{1/2}. 
$$\end{thm}

  \bigskip\par 
Now we turn to upper bounds. We assume that  $d$ is   sub-multiplicative:
\begin{equation}\label{hal0}
 d(nm)\le d(n)d(m) \qq \hbox{provided $(n,m)=1$.}
\end{equation} 
     A
typical example  is for instance function $d_K(n)= \chi\{(n,K)=1\}$. Naturally all multiplicative functions   are
sub-multiplicative, and so is the case of    
 $ d (n)= \l^{\omega(n)}  $,
where  $\l>1$ and $\o(n)= \#\{p: p\mid n\}$.  
\smallskip\par In \cite{W4}   a general upper bound is obtained,   containing and strictly  improving the main results in
\cite{LW1},\cite{LW2}. Further the proof   is entirely based on Gaussian comparison properties,
all suprema of auxiliary Gaussian processes used being   computable exactly.  
 Introduce a basic decomposition. Denote by $P^+(n)$   the largest prime divisor
of $n$.
  Then $$
\{2,\ldots,N\}=\sum_{1\le j\le  \pi(N)}  E_j,  \qq\qq   E_j=\big\{ 2\le n\le N :  P^+(n)=p_j\big\} 
$$
 It is natural to disregard  cells $E_j$ such that $d(n)\equiv0$,  $n\in E_j$. We thus
set
 $\mathcal H_d=\big\{1\le j\le \pi(N):   d_{| E_j}\not\equiv 0\big\}$, $\tau_d=\max \mathcal (H_d)$. 
  The relevant assumption is the following:   
 \begin{equation}\label{pest}p| n\  \Longrightarrow  \ d(n)  \le C\, d({n\over p }),   \qquad  {\rm and}\qquad    d(p^j)\le C_1\l^j,
 \end{equation}
for some   positive  $C,C_1, \l$ with $\l<\sqrt 2$,  any prime
 number $p$, any integers $n,j$.
 
 Clearly, if $C<\sqrt 2$, the second property is implied by the first, although this is not always so as the following example yields.
Fix some prime number $P_1$ as well some reals $1<\l_1<\sqrt 2$, $C_1\ge 1$, and put 
\begin{equation} d  (n)= 
\cases{ C_1\l^j,  & \qq if $P_1^j||\, n$,\cr
 1,  & \qq if $(n,P_1)=1$.
}
\end{equation}
Then $d$ is sub-multiplicative, and satisfies condition (\ref{pest}) with a  constant $C$ which   has to be larger than $  C_1\l$. 
 Now put
 \begin{eqnarray} \label{char}
  D_1(M)&=& \sum_{m=1}^M d(m), \qquad\   \td_1(M)=\max_{1\le m\le M} \frac {D_1(m)}{m}   ,\cr 
   D_2(M) &=& \sum_{m=1}^M d(m)^2,  \qquad  \td_2^2(M)=\max_{1\le m\le M} \frac {D_2(m)}{m} .
 \end{eqnarray}
 
  
\begin{thm} \label{t2} {\gsec (Upper bound)}
 Let   $d$ be a non-negative  sub-multiplicative   function. Assume that condition (\ref{pest}) is realized.
 Let $0\le\s < 1/2$. Then there exists a
constant $C_{\s,d}$ depending on $\s$ and $d$ only, such that   for any integer $N\ge 2$,
 $$
 \E\, \sup_{t \in \R} |\mathcal D (\s+it) |\le  C_{\s,d}  \, \td_2(N)\, B,$$
where $$
B=\cases{
  {N^{1/2-\s}   \tau_d^{1/2} \over (\log N)^{1/2}}
   &,\ {\rm if} \ $\big( {N\log\log N\over \log N}\big)^{1/2}  \le
\tau_d\le \pi(N) ,$  \cr
& \cr
  {N^{3/4-\s} (\log\log N)^{1/4}\over (\log N)^{3/4}}
   &,\ {\rm if} \  $\big( {N \over (\log N)\log\log N}\big)^{1/2}\le
\tau_d \le \big( {N\log\log N\over \log N}\big)^{1/2},$ \cr
& \cr
 N^{1/2-\s}        
          \big({ \tau_d \log\log\tau_d \over  \log\tau_d  }\big)^{1/2}
   &,\ {\rm if} \  $1\le \tau_d \le
\big( {N \over (\log N)\log\log N}\big)^{1/2} .$
} $$
 \end{thm}
  This yields, when combined with  Theorem \ref{lb}, sharp estimates. Consider the following example.

 \smallskip\par
{\it Example 1.}  Take some positive integer
$K$, 
  and  let $d_K(n)= \chi\{(n,K)=1\}$.
Then $d_K$ is   sub-multiplicative and condition (\ref{pest}) is satisfied  with $C=1=\l$. By (\ref{e12}),  this   
defines the remarkable   class of random Dirichlet polynomials,
\begin{equation}
  \DD(s) =\sum_{  (n,K)=1\atop 1\le n\le N } \frac{\e_n }{n^{ s}}  ,
\end{equation}  
  containing   the one of ${\cal E}_\tau$-based Dirichlet
polynomials   considered in \cite{Q3} and \cite{LW1}, where ${\mathcal E}_\tau=\big\{2\le n\le N  : P^+(n)\le p_\tau\}$.  
  Here    $\mathcal H_{d_{K_\tau}}=\sum_{j\le \tau} E_j $.  
We therefore neglect   cells $E_j$, $j>\tau $. Further, we have $\td_1(N) =\td_2(N) \le 1$. As a   consequence of  
Corollary 3 of \cite{W4} and  Theorem \ref{lb}, we have in particular
 
 \begin{thm} Let $0\le\s < 1/2$.  
\smallskip\par\noi   
a) If $\big( {N\log\log N\over \log N}\big)^{1/2}  \le
\tau \le \pi(N) ,$
 $$
 C_{1}(\s) {N^{1/2-\sigma}\tau^{1/2}\over (\log N)^{1/2}}
\le
\E  \sup_{t \in \R} \big|
   \sum_{n\in {\cal E}_\tau}\e_n  n^{ -\s - it}\big|
   \le
   C_{2}(\s) {N^{1/2-\sigma}\tau^{1/2}\over (\log N)^{1/2}}  
  .
$$
   b) If $\big( {N \over (\log N)\log\log N}\big)^{1/2}\le
\tau  \le \big( {N\log\log N\over \log N}\big)^{1/2} $.
 $$
C_1(\sigma)\
{N^{1/2-\s} \tau^{1/2}
\over (\log N)^{1/2}}\
\le
 \E\, \sup_{t \in \R} \big|
   \sum_{n\in {\cal E}_\tau}\e_n  n^{ -\s - it}\big|
  \le C_{2}(\s){N^{3/4-\s} (\log\log N)^{1/4}\over (\log N)^{3/4}}.
$$
  c) If $1\le \tau  \le
\big( {N \over (\log N)\log\log N}\big)^{1/2} .$
 Assume   that $\tau\ge N^{\e}$ for some fixed $0<\e<1/2$.  Then,
$$
C_1(\s,\e)\ {N^{1/2-\s} \tau^{1/2}
\over (\log\tau)^{1/2}}\
\le
\E\, \sup_{t \in \R} \big|
   \sum_{n\in {\cal E}_\tau}\e_n  n^{ -\s - it}\big|
\le
C_{2}(\s)   N^{1/2-\s}        
          \big({ \tau \log\log\tau  \over  \log\tau  }\big)^{1/2}.
$$
 \end{thm}

We notice that the gap is always less that $   (\log\log N)^{1/2}$.
      Theorem \ref{t2}   also applies (see \cite{W4}) to the case $\tau\ll_\e N^{\e}$, as well as  to other classes of examples, for instance  
 \medskip\par 
{\it Example 2.} Consider   multiplicative
functions satisfying the following   condition:
\begin{equation} \label{basic_n}
      \frac{d(p^{a})}{d(p^{a-1 })} \le \l, \qquad  a= 1,2, \ldots
\end{equation}
Clearly (\ref{basic_n})   implies   (\ref{pest}) and further
 $M_d:=\sup_{p} d(p)   <\infty
$,
with $M_d\le \l d(1)$. By theorem 2 of
\cite{HR},   any non-negative multiplicative function $d$ satisfying a Wirsing type condition
$ d(p^m)\le \l_1\l_2^m,
$
for some constants $\l_1>0$ and $0<\l_2<2$ and all prime powers $p^m\le x$, also satisfies
\begin{equation}\label{hr_1}
 {1\over x}\sum_{n\le x} d(n)\le C(\l_1,\l_2) \exp\Big\{\sum_{p\le x} {d(p)-1\over p}\Big\},
\end{equation}
where $C(\l_1,\l_2)$ depends on $ \l_1,\l_2 $ only. This and the fact that
      $d^2$ is
multiplicative and satisfies (\ref{basic_n}) with  $\l^2<2$, yield    that   
  \begin{equation}\label{hr_3}
  \td_1 (N)   \le  C(\l )(\log N)^{M_d}
,\qq \td_2 (N)    \le   C(\l ) (\log N)^{ M_d^2}. 
\end{equation}

\medskip\par 
 

\noi{\gsec Proof of Theorem \ref{t2} \it (Sketch).}
 The proof is    long and technically delicate. We only outline the main steps and will avoid calculation 
details.  
          Let
$M\le N$ and
$0<\s< 1/2$.  
    Fix some integer $\nu$ in $[1,\tau]$ and let  
 $F_\nu  
 = \sum_{1\le j\le \nu} E_j$, $F^\nu=\sum_{  \nu<j\le \tau }E_j $.   
The basic principle of the proof consists of a decomposition  of $Q$    in (\ref{bohr}) into a sum of two   trigonometric
polynomials  
$Q=Q^\e_1+Q^\e_2$, where
$$
Q^\e_1(\uz)  =  \sum_{  n\in F_{\nu }} \e_n d(n) n^{-\s}
e^{2i\pi\langle\ua(n),\uz\rangle}, \
\qq
Q^\e_2(\uz)  =  \sum_{  n\in F^{\nu }}
\e_n  d(n) n^{-\s}e^{2i\pi\langle \ua(n),\uz\rangle}.
$$
 By the contraction principle 
$\E\, \sup_{\uz \in \T^\tau}\big|Q^\e_i(\uz)\big|
\le  C
\E\, \sup_{\uz \in \T^\tau}\big|Q_i(\uz)\big|$, $i=1,2 
$ where $ Q_i $ is the same process as $Q_i^\e $
except that the Rademacher random variables $\e_n$ are replaced by
independent
${\cal N}(0,1)$ random variables $\m_n$. Consequently, both the supremums of $Q_1$ and of $Q_2$ can be estimated, via
  their associated $L^2$-metric.  
  First  evaluate the supremum of
$Q_2$.   We have   
\begin{eqnarray*}
Q_2(\uz) 
&=& \sum_{ \nu <j\le \tau} e^{2i\pi z_j  }\sum_{n\in E_j}
\mu_n  d(n)  n^{-\s}e^{2i\pi \{\sum_{k} a_k({n\over p_j}) z_k \}}.
\end{eqnarray*}
 And so 
\begin{equation}\label{Q2A}\sup_{\uz \in \T^\tau}\big|Q_2(\uz)\big|
\le
4\sup_{\g \in\G}  \big|X (\g)\big|,
\end{equation} 
 where   the
random process $X $ is defined by 
\begin{equation}\label{dec}
X (\g) =\sum_{ \nu <j\le \tau} \a_j
\sum_{n\in E_j}  \mu_n  {d(n)\over  n^{ \s}} \b_{{n\over p_j}},
\qq \g\in \G,
\end{equation}
with
$\g =\big((\a_j)_{\nu <j\le \tau}, (\b_m)_{1\le m\le N/2}\big)
$
and
 $\G=\big\{ \g  :  |\a_j|\vee |\b_m|\le \!\!1,  \nu < \!j\le \tau,
1\! \le m\le N/2\big\} 
 $.  The problem now reduces to estimating the
supremum over $\Gamma$ of the real valued Gaussian process $X$. Plainly
\begin{eqnarray*}
\|X_\g-X_{\g'}\|_2^2
 &\le&
2\!\sum_{ \nu <j\le \tau}\sum_{n\in E_j}  d(n)^2 n^{-2\s} \big[(\a_j -\a'_j)^2+
(\b_{{n\over p_j}}-\b'_{{n\over p_j}})^2\big]  .
\end{eqnarray*}
  Condition (\ref{pest}) and Abel summation yield
\begin{eqnarray*}
& &\|X_\g-X_{\g'}\|_2^2
  \le 
2\!\sum_{ \nu <j\le \tau}\sum_{n\in E_j}  {d(n)^2 \over n^{ 2\s}}\big[(\a_j -\a'_j)^2+
(\b_{{n\over p_j}}-\b'_{{n\over p_j}})^2\big]  \cr &\le &\l^2  \sum_{\nu<j\le\tau} (\a_j -\a'_j)^2 {N^{1-2\s}\td _2^2(N/p_j) \over p_j}+C \l^2\sum_{m \le N/p_\nu} K_m^2
(\b_m-\b'_m)^2,
\end{eqnarray*}
 where $K_m=\sum_{{\nu<j\le \tau\atop  mp_j\le N}}
{ d(m)^2 \over (mp_j)^{2\s}}$\label{Q2D} and
  $\sum_{m\le N/p_\nu} K_m
  \le 
{C N^{1-\s} \td_1({N\over p_\nu})/ \sqrt \nu  \log\nu}   $, by Abel summation. 
  Define a second Gaussian process by putting for all $\g\in \G$
$$
Y(\g) =
\sum_{\nu<j\le\tau} \big({\td_2^2(N/p_j) N^{1-2\s}\over p_j}\big)^{{1/ 2}}\a_j\xi'_j
+
\sum_{m\le N/p_\nu} K_m \ \b_m\xi''_m
 :=
 \ Y'_\g + Y''_\g ,
$$
where $\xi'_i  $, $\xi''_j$  are independent ${\cal N}(0,1)$ random
variables. Thus for some suitable constant
$C$, one has $\|X_\g-X_{\g'}\|_2\le C \|Y_\g-Y_{\g'}\|_2$  for all $\g, \g'\in \G$.
 By  the Slepian  lemma (\cite{W8}, Lemma 10.2.3),  and (\ref{Q2A}) 
\begin{equation}\label{Q2C}
\E\, \sup_{\uz \in \T^\tau}\big|Q_2(\uz)\big|
\le   C  
\E\, \sup_{\g \in \G}\big|Y(\g)\big|.
\end{equation}
As \begin{eqnarray} \label{e24}
 \E\, \sup_{\g\in \G}|Y'(\g)|
 & \le& C\ N^{{1\over 2}-\s}  \td_2(N/p_\nu) \ {\tau^{1/2} \over
(\log\tau)^{1/2}}
 \cr 
 \E\, \sup_{\g\in \G} |Y''(\g)|
&\le & 
{C  N^{1-\s} \td_1(N/p_\nu)\over \nu^{1/2} \log\nu}  ,\end{eqnarray}
   by reporting (\ref{e24})  into (\ref{Q2C}), we get 
\begin{equation}\label{Q2D}
\E\, \sup_{\uz \in \T^\tau}\big|Q_2(\uz)\big|
\le      C\Big(   N^{1/2-\s}  \td_2(N/p_\nu) \ {\tau^{1/2} \over (\log\tau)^{1/2}}+{   N^{1-\s} \td_1(N/p_\nu)\over \nu^{1/2} \log\nu}
\Big) .
\end{equation}

 \medskip

For estimating the supremum of $Q_1 $, we introduce the auxiliary Gaussian process
$$
\Upsilon (\uz) =
\sum_{n\in F_\nu } d(n)  n^{-\s}
\big\{\t_n \cos 2\pi \langle \ua(n),\uz\rangle +
\t_n'\sin 2\pi \langle \ua(n),\uz\rangle \big\}
,\qq  \uz\in \T^{\nu },
$$
where $\t_i$, $\t'_j$  are independent ${\cal N}(0,1)$
random variables. By symmetrization 
\begin{equation}\label{Q1A}
\displaystyle{\E\, \sup_{\uz \in \T^{\nu }}\big|Q_1(\uz)\big|\le
\sqrt{8\pi}
\E\, \sup_{\uz \in  \T^{\nu }}\big|\Upsilon (\uz)\big|}.
\end{equation}
 Plainly $\|\Upsilon(\uz)-\Upsilon(\uz)\big\|_2^2
 \le   4\pi^2\ \sum_{n\in F_\nu }  {d(n)^2\over n^{ 2\s}}
    \big[ \sum_{j=1}^\nu a_j(n) |z_j - z'_j|\big]^2
 .$ 
Now, 
$$\sum_{n\in F_\nu }  {d(n)^2\over n^{ 2\s}}\Big[ \sum_{j=1}^\nu a_j(n) |z_j - z'_j|\Big]^2=\sum_{j =1}^\nu  
      |z_{j } - z'_{j }|^2\sum_{n\in F_\nu }{a_{j }(n)^2d(n)^2\over n^{ 2\s}}$$$$ +\sum_{1\le j_1,j_2\le \nu\atop j_1\not=j_2}    |z_{j_1} - z'_{j_1}|\
|z_{j_2} - z'_{j_2}|\sum_{n\in F_\nu }{a_{j_1}(n)
    a_{j_2}(n)d(n)^2\over n^{ 2\s}}
 :=   S+R.$$

 By sub-multiplicativity, we have for the   rectangle  terms
 \begin{equation}\label{debutR} R     \le     C \sum_{  j_1\not=j_2} |z_{j_1} - z'_{j_1}|  |z_{j_2} - z'_{j_2}|
    \sum_{b_1,b_2=1}^\infty      { b_1d(p_{j_1}^{ b_1})^2 \over  p_{j_1}^{2 b_1\s}   } {b_2 d( p_{j_2}^{ b_2})^2\over  p_{j_2}^{2 b_2\s}
}
  \sum_{k\le {N / p_{j_1}^{ b_1} p_{j_2}^{ b_2}}} {d(k)^2 \over k^{
2\s}  } 
  .
\end{equation}
By using Abel summation, one deduces that
 $$R\le     C_\l \td_2(N)^2 N^{1-2\s} \Big[\sum_{j =1}^\nu {|z_{j } - z'_{j }| \over p_j}
 \Big]^2
 , $$
     As to the   square terms, we have similarly
\begin{eqnarray}\label{debutS}S &\le & \sum_{j =1}^\nu  
      |z_{j } - z'_{j }|^2\sum_{b=1}^\infty {b^2d(p_j^b)^2\over p_j^{ 2b\s}}\sum_{m\le {N\over p_j^b}   }{ d(m)^2\over m^{ 2\s}}\cr& \le
&  C
\td_2(N)^2 N^{1-2\s}\Big[\sum_{j =1}^\nu  
     { |z_{j } - z'_{j }|^2\over p_j} \Big].
 \end{eqnarray}

 Consequently 
 \begin{eqnarray}\label{e27a}
 \|\Upsilon(\uz)-\Upsilon(\uz)\big\|_2&\le& C_\l N^{1/2-\s} \td_2(N)  \max \bigg[  \sum_{j=1}^\nu  {|z_{j} - z'_{j}|\over p_j}
, \Big[\sum_{j =1}^\nu  
     { |z_{j } - z'_{j }|^2\over p_j} \Big]^{1\over 2}\bigg]     \cr 
&\le &    C_\l N^{1/2-\s} \td_2(N)  (\log\log \nu)^{1/2}  \Big( \sum_{j=1}^\nu    {|z_{j} - z'_{j}|^2\over p_j }\Big)^{1/2} 
    , 
\end{eqnarray}
where we used  Cauchy-Schwarz's inequality in the last inequality.
 Let   $g_1,\ldots, g_\nu$ be   independent standard Gaussian r.v.'s. Then
 $U(z):=   \sum_{j=1}^\nu  g_j p_j^{-1/2} z_{j} $ satisfies 
 $\|U(z)-U(z')\|_2=  \big( \sum_{j=1}^\nu    {|z_{j} - z'_{j}|^2\over p_j }\big)^{1/2} . 
 $
 And so 
\begin{equation}\label{slep}\big\|\Upsilon(\uz)-\Upsilon(\uz)\big\|_2\le C_\l N^{1/2-\s} \td_2(N)  (\log\log \nu)^{1/2} \|U(z)-U(z')\|_2.
\end{equation}
 Henceforth
$$\E\,  \sup_{\uz,\uz'\in T^\nu} |\Upsilon (\uz')-\Upsilon(\uz)|\le C_\l N^{1/2-\s} \td_2(N) (\log\log \nu)^{1/2}\E\,  \sup_{\uz,\uz'\in
T^\nu} | U  (\uz')- U (\uz)|
 .  $$ 
 But obviously 
 $ \sup_{\uz \in T^\nu}  |U(z)| =     \sum_{j=1}^\nu  |g_j |p_j^{-1/2}$, and so  
  $\E\,  \sup_{\uz'\in T^\nu} | U  (\uz')- U (\uz)|
  \le C  ({\nu / \log \nu}  )^{1/2}$. By reporting, and since
$\|\Upsilon (\uz)\|_2\le C N^{1/2-\s} \td_2(N) $, for any $\uz\in
\T^{\nu}$, we get 
  \begin{equation} \label{e211}
\E\,  \sup_{\uz'\in  T^\nu} |\Upsilon (\uz')|
\le    C
N^{1/2-\s} \td_2(N) \big({\nu \log\log \nu\over \log \nu} \big)^{1/2}.
\end{equation}

 By substituting in  (\ref{Q1A}) and combining with (\ref{Q2D})    consequently get  
\begin{equation}  \label{coll} \E\, \sup_{t \in \R} |\mathcal D (\s+it) |\le
  C_{\s,\l}  \,N^{1/2-\s}\td_2(N)  \, \bigg[ \big({\nu \log\log \nu\over \log \nu} \big)^{1/2}+
          { \tau^{1/2} \over (\log\tau)^{1/2}}+{   N^{1/2 }  \over \nu^{1/2} \log\nu} \bigg]
 .
\end{equation} 
 The proof is achieved by
  estimating separately the upper bound  in the three cases:
 \par  
{\it  i)}\quad $ ( {N\log\log N\over \log N} )^{1/2}  \le
\tau\le \pi(N).$
 \par    {\it  ii)}\quad $ ( {N \over (\log N)\log\log N} )^{1/2}\le
\tau \le  ( {N\log\log N\over \log N} )^{1/2}.$
 \par     {\it  iii)}\quad   $1\le \tau \le
\big( {N \over (\log N)\log\log N} )^{1/2}.$
  \cqfd
 
\bigskip\par   
\noi  {\sl Acknowledgments.} {\it The author   thanks R. de la Bret\`eche and A. Ivi\'c for   useful  
references.}

 { 
}

\noi {\phh Michel  Weber, \noi  Math\'ematique (IRMA),
Universit\'e Louis-Pasteur et C.N.R.S.,   7  rue Ren\'e Descartes,
67084 Strasbourg Cedex, France.
\par\noindent
E-mail: \  \tt weber@math.u-strasbg.fr}
\end{document}